\documentclass[11pt]{amsart}
\usepackage{amscd,amssymb}
\usepackage[arrow,matrix]{xy}
\newtheorem{theorem}{Theorem}[section]
\newtheorem{lemma}[theorem]{Lemma}
\newtheorem{prop}[theorem]{Proposition}
\newtheorem{corollary}[theorem]{Corollary}
\newtheorem{question}[theorem]{Question}

\theoremstyle{definition}
\newtheorem{definition}[theorem]{Definition}
\newtheorem{example}[theorem]{Example}
\newtheorem{remark}[theorem]{Remark}

\newtheorem*{ack}{Acknowledgment}

\def\ca{{\mathcal A}}
\def\cb{{\mathcal B}}
\def\ck{{\mathcal K}}
\def\h{{\mathfrak h}}
\def\g{{\mathfrak g}}
\def\Z{{\mathbb Z}}
\def\Q{{\mathbb Q}}
\def\C{{\mathbb C}}
\def\F{{\mathbf F}}
\def\k{{\Bbbk}}

\DeclareMathOperator{\gr}{gr}
\DeclareMathOperator{\ext}{Ext}
\DeclareMathOperator{\Hom}{Hom}
\DeclareMathOperator{\tor}{Tor}
\DeclareMathOperator{\fl}{Fl}
\DeclareMathOperator{\Lie}{{Lie}}
\DeclareMathOperator{\Sym}{Sym}
\DeclareMathOperator{\im}{{im}}
\DeclareMathOperator{\rk}{{rk}}

\def\vareps{{\varepsilon}}
\def\set#1{{\left\{#1\right\}}}
\def\abs#1{{\left|#1\right|}}
\def\bfcdot{{\hbox{$\scriptscriptstyle\bullet$}}}

\begin{document}
%\date{January 27, 2005}
%\date{July 19, 2005}

\title[Homotopy Lie algebra of an arrangement]%
{On the homotopy Lie algebra of an arrangement}

\author[G. Denham]{Graham Denham$^1$}
\address{Department of Mathematics, University of Western Ontario,
London, ON  N6A 5B7}
\email{{gdenham@uwo.ca}}
\urladdr{{http://www.math.uwo.ca/\~{}gdenham}}
\thanks{{$^1$}Partially supported by a grant from NSERC of Canada}

\author[A.~I. Suciu]{Alexander~I.~Suciu$^2$}
\address{Department of Mathematics,
Northeastern University,
Boston, MA 02115}
\email{{a.suciu@neu.edu}}
\urladdr{{http://www.math.neu.edu/\~{}suciu}}
\thanks{$^2$Partially supported by NSF grant DMS-0311142}



\begin{abstract}
Let $A$ be a graded-commutative, connected $\k$-algebra 
generated in degree $1$.  The homotopy Lie algebra $\g_A$ 
is defined to be the Lie algebra of primitives of the Yoneda 
algebra, $\ext_{A}(\k,\k)$.   Under certain homological assumptions 
on $A$ and its quadratic closure, we express $\g_A$ as a 
semi-direct product of the well-understood 
holonomy Lie algebra $\h_A$ with a certain $\h_A$-module.  
This allows us to compute the homotopy Lie 
algebra associated to the cohomology ring of the complement 
of a complex hyperplane arrangement, provided some 
combinatorial assumptions are satisfied.   As an application, 
we give examples of hyperplane arrangements whose 
complements have the same Poincar\'e polynomial, 
the same fundamental group, and the same holonomy 
Lie algebra, yet different homotopy Lie algebras.

\end{abstract}

\maketitle

\section{Definitions and statements of results}
\label{sec:intro}

\subsection{Holonomy and homotopy Lie algebras}
\label{sec:} 
Fix a field $\k$ of characteristic $0$.  
Let $A$ be a graded, graded-commuta\-tive algebra over 
$\k$, with graded piece $A_k$, $k\ge 0$.   We will 
assume throughout that $A$ is locally finite, connected, 
and generated in degree $1$.  In other words, $A=T(V)/I$, 
where $V$ is a finite-dimensional $\k$-vector space, 
$T(V)=\bigoplus_{k\ge 0} V^{\otimes k}$ 
is the tensor algebra on $V$, and  $I$ is a 
two-sided ideal, generated in degrees $2$ and higher.   
To such an algebra $A$, one naturally associates two 
graded Lie algebras over $\k$ (see for instance \cite{Av98}, 
\cite{FL}).  

\begin{definition}
\label{def:holonomy}
The {\em holonomy Lie algebra} $\h_{A}$ is the quotient of the 
free Lie algebra on the dual of $A_1$, modulo the ideal generated 
by the image of the transpose of the multiplication map 
$\mu\colon A_1\wedge A_1 \to A_2$:
\begin{equation}
\label{eq:hlie}  
\h_A= \Lie(A_1^*)  \big\slash \text{ideal}\, 
(\im (\mu^*\colon A_2^* \to A_1^* \wedge A_1^*) )
\big. .
\end{equation}
\end{definition}

Note that $\h_A$ depends only on the quadratic closure of $A$:  
if we put $\overline{A}=T(V)/(I_2)$, then $\h_{A}=\h_{\overline{A}}$. 

\begin{definition}
\label{def:homotopy}
The {\em homotopy Lie algebra} $\g_{A}$ is the graded Lie algebra 
of primitive elements in the Yoneda algebra of $A$:
\begin{equation}
\label{eq:glie} 
\g_A=\operatorname{Prim} (\ext _A(\k,\k)).
\end{equation}
\end{definition}

In other words, the universal enveloping algebra of the 
homotopy Lie algebra is the Yoneda algebra:
\begin{equation}
\label{eq:ug}
U(\g_A)=\ext_A(\k,\k).  
\end{equation}

The algebra $U=\ext_A(\k,\k)$ is a bigraded algebra;  let us 
write $U^{pq}$ to denote cohomological degree $p$ and 
polynomial degree $q$.  Then $U^{pq}=0$, unless $-q\geq p$. 
The subalgebra $R=\bigoplus_{p\geq0} U^{p,-p}$ is called the 
linear strand of $U$.  For convenience, we will let 
$U^p_q=U^{p,-p-q}$. The lower index $q$ is called 
the internal degree.  Then $U$ is a graded
$R$-algebra, with $R=U_0$.  Note that 
$U_+=\bigoplus_{q>0} U_{q}$ is an ideal in $U$, 
with $U/U_+ \cong R$.

The relationship between the holonomy and homotopy Lie 
algebras of $A$ is provided by the following well-known 
result of L\"ofwall. 

\begin{lemma}[L\"ofwall \cite{lof86}]
\label{lem:lof}
The universal enveloping algebra of the holonomy Lie 
algebra, $U(\h_A)$, equals the linear strand, 
$R=\bigoplus_{p\geq0} U^p_0$, of the Yoneda 
algebra $U=U(\g_A)$.
\end{lemma}

Particularly simple is the case when $A$ is a Koszul algebra.  
By definition, this means the homotopy Lie algebra $\g_A$ 
coincides with the holonomy Lie algebra $\h_A$, i.e., $U=R$. 
Alternatively, $A$ is quadratic (i.e., $A=\overline{A}$), and its 
quadratic dual, $A^!=T(V)/(I_2^{\perp})$, coincides with the 
Yoneda algebra: $A^{!}=U$.  For an expository account of 
Koszul algebras, see \cite{frob99}.

As a simple (yet basic) example, take $E=\bigwedge V$, 
the exterior algebra on $V$.  Then $E$ is 
Koszul, and its quadratic dual is $E^{!}=\Sym(V^*)$, the 
symmetric algebra on the dual vector space.  Moreover, 
$\g_A=\h_A$ is the abelian Lie algebra on $V$. 

\subsection{Main result}
\label{subsec:main result}
The computation of the homotopy Lie algebra of 
a given algebra $A$ is, in general, a very hard problem. 
Our goal here is to determine $\g_A$ under certain 
homological hypothesis.   First, we need one more definition.

Let $B=\overline{A}$ be the quadratic closure of $A$.   
View  $J=\ker (B\twoheadrightarrow A)$ as a graded left 
module over $B$. 

\begin{definition}
\label{def:homotopy module}
The {\em homotopy module} of a graded algebra $A$ is 
\begin{equation}
\label{eq:M module}
M_A=\ext_B(J,\k), 
\end{equation}
viewed as a bigraded left module over the ring 
$R=U(\h_A)=\ext_B(\k,\k)$ via the Yoneda product.  
\end{definition}

\begin{theorem}
\label{th:main}
Let $A$ be a graded algebra over a field $\k$, with quadratic 
closure $B=\overline{A}$, and homotopy module $M=M_A$.  
Assume $B$ is a Koszul algebra, and there exists an integer 
$\ell$ such that $M_q=0$ unless $\ell\leq q\leq\ell+1$.  
Then, as graded Hopf algebras, 
\begin{equation}
\label{eq:UT}
U(\g_A)\cong T(M_A[-2])\, \widehat{\otimes}_\k\, U(\h_A).  
\end{equation} 
\end{theorem}

Here $M[q]$ is the graded $R$-module  
with $M[q]^r=M^{q+r}$, while 
$T(M[-2]) \widehat{\otimes}_\k R$ is the ``twisted" 
tensor product of algebras, with underlying vector space 
$T(M[-2]) \otimes_\k R$ and multiplication 
$(m\otimes r)\cdot (n\otimes s) = 
(-1)^{\abs{r}\abs{n}} ((m\otimes n) \otimes rs + 
(m\otimes nr) \otimes s)$.  
%where the action of $R$ 
%on $T(M[-2])$ is induced from the 
%$R$-module structure of $M[-2]$.

Taking the Lie algebras of primitive elements in the 
respective Hopf algebras, we obtain the following.

\begin{corollary}
\label{cor:lie}
Under the above hypothesis, the homotopy Lie algebra of $A$ 
splits as a semi-direct product of the holonomy Lie algebra 
with the free Lie algebra on the (shifted) homotopy module,
\begin{equation}
\label{eq:glie pres}
\g_A\cong \Lie(M_A[-2])\rtimes \h_A,
\end{equation}
where the action of $\h$ on $\Lie(M)$ is given by 
$[m,h]=-hm$ for $h\in\h$ and $m\in M$.
\end{corollary}

As pointed out to us by S.~Iyengar, Theorem \ref{th:main} 
implies (under our hypothesis) that the projection map 
$U(\g_A)\to U(\h_A)$ is a Golod homomorphism.  
Therefore, the semi-direct product structure of $\g_A$ also 
follows from results of Avramov \cite{Av78}, \cite{Av86}.

\subsection{Hyperplane arrangements}
\label{subsec:arr}
Let $\ca=\{H_1,\dots, H_n\}$ be an arrangement of 
hyperplanes in $\C^{\ell}$, with intersection lattice $L(\ca)$ 
and complement $X(\ca)$. 
The cohomology ring $A=H^\bfcdot(X(\ca),\k)$ admits a 
combinatorial description (in terms of $L(\ca)$), due to 
Orlik and Solomon: 
\begin{equation}
\label{eq:os}
A=E/I,
\end{equation}
 where $E$ is the exterior algebra over $\k$, on generators 
 $e_1, \dots , e_n$ in degree $1$, and $I$ is the ideal 
 generated by all elements of the form 
 $\sum_{q=1}^{r}(-1)^{q-1}e_{i_1} \cdots \widehat{e_{i_q}} 
 \cdots e_{i_r}$ for which 
 $\rk (H_{i_1}\cap \cdots \cap H_{i_r}) < r$;  see \cite{ot}. 

The holonomy Lie algebra of the Orlik-Solomon algebra 
also admits an explicit presentation, this time solely in terms 
of $L_{\le 2}(\ca)$.  Identify $\Lie(A_1^*)$ with the free Lie 
algebra over $\k$, on generators $x_H=e_H^*$,  $H\in \ca$.  
Then:
\begin{equation}
\label{eq:hlie arr}  
\h_A= \Lie (A_1^*) \Big\slash \text{ideal}\, \Big\{
\Big[x_H , \sum_{H'\in \ca \colon H'\supset F} x_{H'}\Big] \:\big| \:
F \in L_2(\ca) \text{ and } F\subset H \big. \Big\}  \Big. .
\end{equation}

As we shall see in Section~\ref{sec:lie}, the homotopy 
Lie algebra $\g_A$ also admits a finite presentation, for 
a certain class of hypersolvable arrangements, to be defined below. 

\begin{question}
\label{question:fin pres}
Do there exist arrangements for which $\g_A$ is not finitely
presented?  For which the (bigraded) Hilbert series of $U(\g_A)$ is
not a rational function?
\end{question}

\subsection{Hypersolvable arrangements}
\label{subsec:hypr arr}
An arrangement $\ca$ is called {\em supersolvable} if its intersection 
lattice admits a maximal modular chain.  The OS algebra of a 
supersolvable arrangement has a quadratic Gr\"obner basis, and  
thus, it is a Koszul algebra (this result, implicit in Bj\"orner and 
Ziegler \cite{BZ}, was proven in Shelton and Yuzvinsky \cite{shyu97}). 

An arrangement $\ca$ is called {\em hypersolvable} if it has the 
same intersection lattice up to rank $2$ as that of a supersolvable 
arrangement.  This ``supersolvable deformation," $\cb$, 
is uniquely defined, and has the property that the two  
complements have isomorphic fundamental groups;  
see Jambu and Papadima \cite{japa98, japa02}.   
Let  $A=H^\bfcdot(X(\ca),\k)$ and $B=H^\bfcdot(X(\cb),\k)$ 
be the respective OS algebras. 
It is readily seen that $B=\overline{A}$; thus, $A$ and $B$ share 
the same holonomy Lie algebra: $\h=\h_{A}=\h_{B}$.  Furthermore, 
since $B$ is Koszul, we have $\g_B=\h$. 

The hypothesis of Theorem~\ref{th:main} holds in two nice 
situations, which can be checked combinatorially; 
for precise definitions, see \S\ref{subsec:sing range} 
and \S\ref{subsec:hyper slice}, respectively. 

\begin{theorem}
\label{thm:main-arr}
Let $\ca$ be an arrangement, and let $A$ be its Orlik-Solomon 
algebra.  Suppose either 
\begin{enumerate}
\item \label{i}  
$\ca$ is hypersolvable, and its singular range 
has length $0$ or $1$; or 
\item \label{ii}  
$\ca$  is obtained by fibred extensions 
of a generic slice of a supersolvable arrangement.  
\end{enumerate}
Then $\g_A\cong \Lie(M_A[-2])\rtimes \h_A$.
\end{theorem}

An explicit finite presentation for $\g_A$ is given in 
Theorem \ref{thm:g pres}, in the case when $\ca$ is 
a generic slice of a supersolvable arrangement.  
The Eisenbud-Popescu-Yuzvinsky resolution \cite{epy03} 
permits us to compute the Hilbert series of $M_A$ (and hence, 
that of $\g_A$) in the case when $\ca$ is a $2$-generic 
slice of a Boolean arrangement.

Theorem~\ref{thm:main-arr} allows us to distinguish between 
hyperplane arrangements whose holonomy Lie algebras are 
isomorphic.  In Example \ref{ex:2gen-7}, we exhibit a pair of 
$2$-generic, $4$-dimensional sections of the Boolean 
arrangement in $\C^7$; the two arrangements have the same 
fundamental group, the same Poincar\'e polynomial, and the 
same holonomy Lie algebra, yet different homotopy Lie algebras.

In Section \ref{sect:top}, we provide some topological interpretations. 
As noted in \cite{CCX}, \cite{ps04}, the holonomy Lie algebra of a 
supersolvable arrangement equals, up to a rescaling 
factor, the topological homotopy Lie algebra of the 
corresponding ``redundant" subspace arrangement.  
We extend this result, and relate the homotopy Lie 
algebra of an arbitrary hyperplane arrangement to 
the topological homotopy Lie algebras of the redundant 
subspace arrangement.  As a consequence, we find a 
pair of codimension-$2$ subspace arrangements in 
$\C^8$, whose complements are simply-connected 
and have the same homology groups, yet distinct higher 
homotopy groups.

\section{Some homological algebra}
\label{sec:homological algebra}

\subsection{The homotopy module}
\label{subsec:jm}
Let $A$ be graded, graded-commutative, connected, locally 
finite algebra.  Assume $A$ is generated in degree $1$, 
and its quadratic closure, $B=\overline{A}$ is a Koszul 
algebra.   Let $E$ be the exterior algebra on $A_1=B_1$.  
Let $I$ and $J$ be, respectively, the kernels of the 
natural surjections $E\twoheadrightarrow B$ and 
$B\twoheadrightarrow A$, giving the exact sequences 
\begin{align}
\label{eq:one}
&\xymatrix{0\ar[r]& I \ar[r]&  E \ar[r]&  B \ar[r]&  0},\\
\label{eq:two}
&\xymatrix{0 \ar[r]&  J \ar[r]&  B \ar[r]&  A \ar[r]&  0}.
\end{align}

In what follows, we will record some homological properties of 
the ring $A$, viewed as a $B$-module.  Recall if $N$ is a 
$B$-module, the Yoneda product gives $\ext_B(N,\k)$ the 
structure of a left module over the ring $R=U(\h_A)=\ext_B(\k,\k)$.  
An object of primary interest for us will be the {\em homotopy 
module} of $A$, 
\begin{equation}
\label{eq:defM}
M=M_A=\ext_B(J,\k).   
\end{equation}
This bigraded $R$-module will play a crucial 
role in the determination of the homotopy Lie algebra $\g_A$. 

Our grading conventions shall be as follows.  
Suppose $V$ and $W$ are $\Z$-graded $\k$-vector spaces.  Then 
$f\in\Hom_\k(V,W)$ has degree $r$ if $f \colon V^q\to W^{q+r}$ 
for all $q$.  For any $\Z$-graded $\k$-vector space $V$,
we shall let $V^*$ denote the graded $\k$-dual of $V$.  
In particular, then, $(V^*)^q=\Hom_\k(V^{-q},\k)$.
If $V$ has finite $\k$-dimension in each graded piece, 
then $(V^*)^*\cong V$.

We shall treat all boundary maps in chain complexes as having 
polynomial degree $0$ and homological degree $+1$.  Then, in 
particular, chain complexes will be regarded as cochain complexes 
in negative degree.  We shall indicate 
shifts of polynomial grading by defining $V(q)^r=V^{q+r}$, and shifts 
of homological grading by writing $V[q]$ analogously.

Following these conventions, $M^{pq}=\ext_B^p(J,\k)^q$ is nonzero only
for $q\leq -p$.  Then, taking $M^p_q=M^{p,-p-q}$ (the internal grading), 
we have $M^p_q\neq0$ only for $q\geq0$.  The grading is such that, 
for each fixed $q$,
the action of $R$ on $M$ satisfies $R^r\otimes M^p_q\rightarrow M^{r+p}_q$.

\begin{lemma}
\label{lem:ExtAk}
$\ext_B(A,\k)\cong \k\oplus M[-1]$ as graded $R$-modules.
\end{lemma}
\begin{proof}
Consider the long exact sequence for $\ext_B(-,\k)$ applied to
\eqref{eq:two}:
\begin{equation}
\label{eq:les}
\xymatrixcolsep{16pt}
\xymatrix{
\cdots\ar[r] & \ext_B^{q-1}(J,k)\ar[r] & \ext_B^q(A,k)\ar[r] &
\ext_B^q(B,k)\ar[r] & \cdots
}
\end{equation}
Since $\ext^0_B(A,\k)\cong\ext_B^0(B,\k)\cong\k$ and $\ext^q_B(B,\k)=0$ for
all $q>0$, the map $\ext_B(B,\k)\to\ext_B(J,\k)$ is zero.  So the long
exact sequence breaks into short exact sequences which,  using 
\eqref{eq:defM}, we will write as a single short exact sequence 
of graded $R$-modules,
\begin{equation}
\label{eq:four}
\xymatrixcolsep{16pt}
\xymatrix{
0\ar[r] & M[-1]\ar[r] & \ext_B(A,\k)\ar[r] & \k\ar[r] & 0
}.
\end{equation}
For each $q$, one of the two maps is zero and the other is an isomorphism,
so the short exact sequence splits.
\end{proof}

\subsection{Injective resolutions}
\label{subsec:injres}
For any $E$-module $N$, let
\begin{equation}
\label{eq:xcirc}
N^{\circ}=\set{a\in E:ax=0\hbox{~for all $x\in N$}},
\end{equation}
the annihilator of $N$ in $E$.  Later on, 
we require explicit, injective resolutions.

\begin{lemma}
\label{lem:injres}
Suppose the ring $B=E/I$ is an arbitrary quotient of a 
finitely-generated exterior algebra $E$.  If 
\begin{equation}\label{eq:fres}
\xymatrixcolsep{16pt}
\xymatrix{
0 & \k\ar[l] & B\otimes_\k F^0\ar[l] & B\otimes_k F^1\ar[l] & \cdots\ar[l]
}
\end{equation}
is a minimal, free resolution of $\k$ over $B$, then 
\begin{equation}\label{eq:ires}
\xymatrixcolsep{16pt}
\xymatrix{
0\ar[r] & \k\ar[r] & B^*\otimes_\k (F^0)^*\ar[r] 
& B^*\otimes_k (F^1)^*\ar[r]
&\cdots}
\end{equation}
is an injective resolution of $\k$ over $B$.
\end{lemma}
\begin{proof}
The resolution \eqref{eq:fres} is an acyclic complex of $E$-modules,
so its vector space dual \eqref{eq:ires} is an acyclic complex as well, 
since each $F^i$ has finite $\k$-dimension.

Now $B^*\cong I^\circ(n)$ as $E$-modules, via the determinantal 
pairing in $E$.  On the other hand, $E$ is injective as a module 
over itself, so $I^\circ$ is injective as an $E$-module; 
see \cite[Prop. 2.27]{shvabook}.  Since each $F^i$ has finite 
$\k$-dimension, each $B^*\otimes_\k (F^i)^*$ is injective.
\end{proof}

\begin{lemma}
\label{lem:injres2}
Let $A$ and $B$ two algebras, with $\overline{A}=B$ 
Koszul.  Write $B=E/I$, $A=B/J$,  
$\h=\h_A=\h_B$, and $R=U(\h)$.  Then:
\begin{enumerate}
\item \label{claim1}
The complex
\begin{equation}
\label{eq:complex}
\xymatrixcolsep{16pt}
\xymatrix{
0\ar[r] & \k\ar[r] & I^\circ(n)\otimes_\k R^0\ar[r] & 
I^\circ(n)\otimes_\k R^1\ar[r] & \cdots
}
\end{equation}
is an injective resolution of $\k$ over $B$, with boundary map 
described below.
\item  \label{claim2}
$\ext^q_B(A,\k)\cong H^q(J^\circ(n)\otimes_\k R^\bfcdot)$, 
for all $q\geq0$.
\end{enumerate}
\end{lemma}
\begin{proof}
The Koszul complex $\ck^*=B\otimes_\k R^*$ is a free $B$-module 
resolution of $\k$, so it is also an acyclic complex of $E$-modules, 
with boundary map induced from
\begin{equation}
\label{eq:partial star}
\partial^* \colon 1\otimes x_i^*\mapsto e_i\otimes 1.
\end{equation}
Then $\Hom_\k(B\otimes_\k R^*,\k)=B^*\otimes_\k R$ is an injective
resolution, by the previous Lemma.

To establish \eqref{claim2}, it suffices to note that
$\Hom_B(A,I^\circ)\cong J^\circ$.
\end{proof}

\section{Proof of the main result}
\label{sect:main}

Our approach to the proof of Theorem~\ref{th:main} 
is to construct a spectral sequence comparing the
minimal resolution and the Koszul complex of $A$.  
We show the spectral 
sequence collapses at $E_2$ under suitable hypotheses
in Proposition~\ref{prop:collapses}, though not in general 
(Example~\ref{ex:lived2}).  This collapsing is enough to prove the 
theorem, via Proposition~\ref{prop:ses}.

\subsection{A spectral sequence}
\label{subsec:spec seq}
Using the previous notation,
$A\otimes_\k U^*\to\k\to0$ is a minimal free resolution of $\k$ over
$A$.  It is filtered by degree, and the linear strand is 
$A\otimes_\k R^*$.
That is, there is a short exact sequence of chain complexes
\begin{equation}
\label{eq: ru}
\xymatrix{
0 \ar[r]& A\otimes_\k R^*  \ar^{1\otimes\vareps^*}[r]& 
A\otimes_\k U^*  \ar[r]& A\otimes_\k U^*_+    \ar[r]&  0
}.
\end{equation}

Now $B\otimes_\k R^*$ is a free resolution of $\k$ over $B$, 
since $B$ is Koszul.  Using Lemma~\ref{lem:ExtAk}, 
we find that the homology of the linear strand 
(Koszul complex) is
\begin{align}
\label{eq: harb}
H_{\bfcdot}(A\otimes R^*)&\cong\tor^B(A,\k) \notag \\
&\cong\ext_B(A,\k)^*\\
&\cong\k\oplus M[-1]^*.  \notag
\end{align}
The long exact sequence in homology then reveals that
\begin{equation}
\label{eq: haum}
H_{\bfcdot}(A\otimes_\k U_+^*)\cong M[-2]^*
\end{equation}
as $A$-modules.  Recall that $A$ acts trivially on $M$ 
(and hence on $M[-2]^*$), so 
\begin{equation}
\label{eq:AU1}
\Hom_A(H_{\bfcdot}(A\otimes_\k U_+^*),\k)\cong M[-2].
\end{equation}

On the other hand, since our complex is a quotient 
of a minimal resolution,
\begin{equation}
\label{eq:AU2}
H_{\bfcdot}(\Hom_A(A\otimes_\k U_+^*,\k))\cong U_+.
\end{equation}
Comparing the two gives a universal coefficients 
spectral sequence of the form
\begin{equation}
\label{eq:E2}
E_2^{pq}=\ext_A^p((M[-2]^*)^q,\k) 
\cong  M[-2]^q\otimes_\k U^p\ \Longrightarrow \ 
U_+^{p+q}.
\end{equation}
The spectral sequence is used as follows.
\begin{prop}\label{prop:ses}
If $E_\infty=E_2$ in the spectral sequence \eqref{eq:E2}, then
\[
\xymatrix{
0\ar[r]& U\otimes_\k M[-2] \ar^(.7){\phi}[r] & U 
\ar^{\vareps}[r]& R \ar[r]&0
}
\]
is exact, and the conclusion of Theorem~\ref{th:main} holds.
\end{prop}
\begin{proof}
If $E_\infty=E_2$, then $U\otimes M[-2]\cong U_+$ as a (left) 
$U$-module.  Now $U_+=\ker\vareps$, giving the short exact 
sequence.  Since $\h$ is a Lie subalgebra of $\g$, $R=U(\h)$ 
is a Hopf subalgebra of $U=U(\g)$, so the sequence splits.  
The isomorphism of Theorem~\ref{th:main} can then
be obtained by induction.
\end{proof}

\subsection{Collapsing conditions}
\label{subsec:collapse}
In order to show that the higher differentials in the spectral sequence
\eqref{eq:E2} vanish, we use a degree argument that begins 
by considering the $E_0$ term.  Since
\begin{equation}
\label{eq:injres}
\xymatrix{
0  \ar[r]&  \k  \ar[r]&  A^*\otimes_\k U^0 \ar[r]&  
A^*\otimes_\k U^1 \ar[r]&  \cdots
}
\end{equation}
is an injective resolution of $\k$ over $A$, (Lemma~\ref{lem:injres})
we consider the double complex
\begin{align}
\label{eq:doublecomplex}
C^{pq}&=\Hom_A(A\otimes_\k (U^q)^*_+,A^*\otimes_\k U^p)\\
&\cong  U^q_+\otimes_\k A^*\otimes_\k U^p,\notag
\end{align}
with induced boundary maps $\partial_h$ and $\partial_v$.  Then our
spectral sequence \eqref{eq:E2} is obtained by filtering $C^{\bfcdot\bfcdot}$
by columns.  Checking the grading, we see
\begin{equation}
\label{eq:dv2}
\partial_v\colon U^q_+\otimes_\k (A^*)^s\otimes_\k U^p
\rightarrow U^{q+1}_+\otimes_\k (A^*)^{s+1}\otimes_\k U^p
\end{equation}
and 
\begin{equation}
\label{eq:dh2}
\partial_h\colon U^q_+\otimes_\k (A^*)^s\otimes_\k U^p
\rightarrow U^q_+\otimes_\k (A^*)^{s+1}\otimes_\k U^{p+1}.
\end{equation}
By looking at $E_2$ and $\partial_v$, we see that we must have $E_1=E_2$.

We first consider the case where the ideal $J$ has a (shifted) linear
resolution.
\begin{prop}
\label{prop:collapses}
Suppose $\ca$ is a hypersolvable arrangement 
for which $M^p_q=0$ unless $q=\ell$, for 
some fixed $\ell$. Then $E_2=E_\infty$.
\end{prop}

\begin{proof}
In this case, $M[-2]^q_r=0$ unless $r=\ell-2$.
Then $H^q(C^{p\bfcdot},\partial_v)_r=0$ unless $r=\ell-2$.

First we note that $(U_+)^q_t=0$ unless $t\geq\ell-2$.  
This can be seen from the fact that $U_+$ is a graded 
subquotient of $M[-2]\otimes_\k U$, from \eqref{eq:E2}: 
the support of $M[-2]$ is described above, and
$U^p_q=0$ unless $q\geq 0$.

Regard $A^*$ as a chain complex concentrated in homological degree $0$.
Then observe that the internal degree of a nontrivial cocycle 
representative in $(U_+)^q_t\otimes_\k (A^*)_s$
is $s+t =\ell-2$, by the first observation above. It follows $s\le 0$ 
from the inequality above.  However, $(A^*)_s=0$ 
unless $0\leq s\leq\ell$, so the representative of a nonzero,
homogeneous $E_2$-cocycle in $E_0$ must have $s=0$.

Now suppose $x\in E_2^{pq}$ is such a cocycle, with representative
$\widetilde{x}$ in $C^{pq}$.  By the above, $\widetilde{x}\in U^q_+
\otimes_\k (A^*)_0$.  Then $\partial_h(\widetilde{x})=0$ in $C^{p+1,q}$
by \eqref{eq:dh2}.  This means $d_2(x)=0$, and similarly for higher
differentials.
\end{proof}

\begin{proof}[Proof of Theorem~\ref{th:main}]
In view of Proposition~\ref{prop:ses}, it
remains only to show the spectral sequence collapses when $M^p_q=0$
unless $0\leq \ell-q \leq1$ for some $\ell$.  In this case, let 
$N=M_\ell $ denote the  $R$-submodule of $M$ of internal 
degree $\ell$.

By the same reasoning as in the proof of Proposition~\ref{prop:collapses},
$N[-2]\otimes U\subseteq \ker d_k$ for $k\geq2$.  Now 
$N[-2]\cong N[-2]\otimes U^0$ is a submodule of the $p=0$ 
column of $E_2$.  Since it is (trivially) not in the image of any 
nonzero differentials, $N[-2]$ is an $R$-submodule of $U$.  

Let $K$ denote the Hopf subalgebra of $U$ generated by 
$R$ and $N[-2]$.  By \cite[Theorem 4.4]{mimo65}, $U$ is 
a free $K$-algebra.  It follows that
$K\cong T(N[-2])\otimes_\k R$.  In the notation of the 
previous proposition, any nontrivial differential $d_k$ with 
$k\geq2$ would lift in $E_0$ to a map 
$U_+\otimes (A^*)_1\otimes U\rightarrow
U_+\otimes (A^*)_0\otimes U$.  We have shown that the targets of 
these maps are unchanged between $E_2$ and $E_\infty$, so it 
follows that the maps themselves must also all be zero.  
\end{proof}

\subsection{A non-collapsing spectral sequence}
\label{subsec:non collapsing}

Calculations with the Macaulay~2 package \cite{GS-m2} 
show that the hypotheses of Theorem~\ref{th:main} 
cannot in general be relaxed: differentials in the spectral 
sequence \eqref{eq:E2} may not be zero.

\begin{example}
\label{ex:lived2}
Consider the arrangement defined by the polynomial 
\[
Q=xyz(x-w)(y-w)(z-w)(x-u)(y-u).
\]  
Let $A$ be the Orlik-Solomon algebra, and 
$M=M_A$ its homotopy module.  It is readily 
seen that $M_q\ne 0$ for $q=3,4,5$.  
An Euler characteristic calculation shows that the 
spectral sequence \eqref{eq:E2} must have a 
nonzero differential 
\[
d_2^{04}\colon M[-2]^4_6\otimes_\k U^0
\to M[-2]^3_5\otimes_\k U^2.
\]
It follows that the Hopf algebra $U(\g_A)$ will not have 
the structure we find in Theorem~\ref{th:main}.
\end{example}

\section{Hypersolvable arrangements}
\label{sec:hyper}

In this section, we apply our main result to certain 
classes of hypersolvable arrangements.  

\subsection{Solvable extensions}
\label{subsec:solv}
We start by reviewing in more detail the notion of a hypersolvable 
arrangement, introduced by Jambu and Papadima in \cite{japa98}.  
Roughly, a hypersolvable arrangement is a linear projection of 
a supersolvable arrangement that preserves intersections 
through codimension two.

\begin{definition}[\cite{japa98}]
An arrangement $\ca$ is hypersolvable if there exist
subarrangements 
$\{0\}=\ca_1\subset \ca_2\subset \cdots\subset \ca_m=\ca$, 
so that each inclusion $\ca_i\subset\ca_{i+1}$ is {\em solvable}. 
In turn, an inclusion of hyperplane arrangements $\ca\subset\cb$ 
is called a solvable extension if:
\begin{enumerate}
\item There are no hyperplanes $H\in\cb\setminus\ca$ and 
$H',H''\in\ca$ with $H'\neq H''$ and $\rk(H\cap H'\cap H'')=2$;
\item 
For any $H,H'\in\cb\setminus\ca$, there is exactly one $H''\in\ca$ with
$\rk(H\cap H'\cap H'')=2$, denoted by $f(H,H')$;
\item 
For any $H,H',H''\in\cb\setminus\ca$, one has 
$\rk(f(H,H')\cap f(H,H'')\cap f(H',H''))\leq 2$.
\end{enumerate}
\end{definition}
It turns out that if $\ca$ is hypersolvable with a sequence of solvable
extensions as above, then for all $i$, the rank of 
$\ca_i$ and $\ca_{i+1}$ differ by at most one.  If the ranks are equal,
the extension is said to be singular; otherwise, the extension is 
nonsingular (or fibred, in the sense of Falk and Randell, \cite{falkrand85}).

If $s$ denotes the number of singular extensions, then, $\rk\ca=m-s$.
Jambu and Papadima show in \cite{japa02} that one can replace
the singular extensions by nonsingular ones in order to construct a 
supersolvable arrangement $\cb$ of rank $m$ that 
projects onto $\ca$, preserving the intersection lattice through 
rank $2$.  That is,

\begin{theorem}
\label{th:slice}
An arrangement $\ca$ is hypersolvable iff there exists a 
supersolvable arrangement $\cb$ and a linear subspace 
$W$ for which $\ca=\cb\cap W$ 
and $L(\ca)_{\leq2}\cong L(\cb)_{\leq2}$.
\end{theorem}
\begin{proof}
The implication ``$\Rightarrow$'' is Theorem~2.4 of \cite{japa02}.  
The converse, due to Jambu (private communication), runs as follows.  
Suppose $\cb$ is supersolvable and there exists a subspace $W$ 
as above.  By definition, $\cb$ has a maximal modular chain 
$F_1<F_2<\cdots<F_m$.  Putting $\cb_i=\cb_{X_i}$ gives a 
sequence of solvable extensions for $\cb$, all fibred.  For 
$1\leq i\leq m$, let $\ca_i=\cb_i\cap W$.  Since collinearity 
relations are preserved, each $\ca_i\subset \ca_{i+1}$ 
is also a solvable extension, so $\ca$ is hypersolvable.
\end{proof}

We remark that, in the above proof,  $\ca_i\subset\ca_{i+1}$ 
is a singular extension if and only if $F_i\cap W=F_{i+1}\cap W$.
The arrangement $\cb$ in called the {\em supersolvable 
deformation} of $\ca$.  For example, any arrangement $\ca$ for 
which no three hyperplanes intersect in codimension three is 
hypersolvable, and its supersolvable deformation is the Boolean 
arrangement in $\C^n$, where $n=\lvert \ca \rvert$. 

\begin{lemma}
\label{lem:free}
Suppose $\ca'\subset\ca$ is a fibred extension.  The projection 
$p\colon X(\ca)\rightarrow X(\ca')$ induces an inclusion 
$A'\hookrightarrow A$ of the respective Orlik-Solomon 
algebras which makes $A$ into a free $A'$-module of 
rank $k=\lvert \ca \setminus \ca'\rvert$. 
\end{lemma}
\begin{proof}
The projection $p\colon X\rightarrow X'$ is a bundle map, 
with fiber $\C\setminus \{\text{$k$ points}\}$.  As noted by 
Falk and Randell \cite{falkrand85}, this bundle admits a section, 
and thus the Serre spectral sequence collapses at the $E_2$ term. 
Hence, $H^\bfcdot(X)\cong H^\bfcdot(X') \otimes 
H^\bfcdot(\vee^k S^1)$.  The result follows.
\end{proof}

\subsection{Singular range} 
\label{subsec:sing range}
We now give some easy to check combinatorial 
conditions insuring that a hypersolvable arrangement 
satisfies the hypothesis of Theorem \ref{th:main}.  
We start by attaching a pair of relevant integers to a 
hypersolvable arrangement.

\begin{definition}
\label{def:srange}
Suppose $\ca$ is hypersolvable with supersolvable deformation 
$\cb$, and $\ca\neq\cb$.  Let $c$ be the least integer
for which $L(\ca)_{\leq c}\not\cong L(\cb)_{\leq c}$.  
Since $\ca\neq\cb$, there is a largest integer $i$ for 
which the extension $\ca_i\subset\ca_{i+1}$ is singular.  
Let $d$ the rank of these two arrangements.  We will 
call the pair $(c,d)$ the {\em singular range} of the 
arrangement $\ca$, and $\abs{d-c}$ the length of this range.
\end{definition}

\begin{lemma}
\label{lem:cd}
If $\ca$ is hypersolvable with singular range $(c,d)$, then $3\leq c\leq d$.
\end{lemma}

\begin{proof}
The inequality $c\geq3$ follows from Theorem~\ref{th:slice}.
Suppose $d<c$; then $L(\ca)_{\leq d}\cong L(\cb)_{\leq d}$.  It follows
that $L(\ca_{d+1})_{\leq d}\cong L(\cb_{d+1})_{\leq d}$, whence 
$\ca_{d+1}=\cb_{d+1}$ since the arrangements are central.  
Since $d$ is greater than or equal to the index of the last singular
extension, however, $\ca_i=\cb_i$ for $d+1\leq i\leq m$, 
so $\ca=\cb$, a contradiction.
\end{proof}

 Let $A=H^\bfcdot(X(\ca),\k)$ 
and $B=H^\bfcdot(X(\cb),\k)$ be the respective Orlik-Solomon algebras.  
Since $L(\ca)_{\leq2}\cong L(\cb)_{\leq2}$, and since the Orlik-Solomon
algebra of a supersolvable arrangement is quadratic, the algebra 
$B=E/I$ is the quadratic closure of $A$. Let 
$J=\ker(B\twoheadrightarrow A$), and let $M=\ext_B(J,\k)$, 
viewed as a module over $R=\ext_B(\k,\k)$.
Since $\cb$ is supersolvable, the algebra $B$ is Koszul 
(see \cite{shyu97});   thus,  $R=B^!$.  

\begin{lemma}
\label{lem:CMreg}
If $\ca$ is a hypersolvable arrangement with singular range $(c,d)$, 
then $M^p_q=0$ unless $p\geq 0$ and $c\leq q\leq d$.
\end{lemma}

\begin{proof}
The ideal $J$ has a minimal, (infinite) free resolution 
over $B$ of the form
\begin{equation}
\label{eq:jresolution}
\xymatrix{
0& J\ar[l]&B\otimes_\k (M^{0,-})^*\ar[l] & B\otimes_\k (M^{1,-})^*\ar[l] 
& \cdots\ar[l].
}
\end{equation}
Recall that $J$ is generated by Orlik-Solomon relations.  By 
Definition~\ref{def:srange}, the least degree of a generator of $J$ is
$c$, so $M^0_c\neq0$ and $M^0_q=0$ for $q<c$.  Thus 
$M^p_q=0$ for $q<c$, establishing the first inequality.

To show $M^p_q=0$ for $q>d$, too, let $i$ be the largest index of a 
singular extension $\ca_i\subset\ca_{i+1}$.
let $B_{i+1}=H^\bfcdot(X(\cb_{i+1}),\k)$ and 
$A_{i+1}=H^\bfcdot(X(\ca_{i+1}),\k)$, and let $B'_{i+1}=
H^\bfcdot(X(\cb'_{i+1}),\k)$ be the cohomology ring 
of the projectivization (decone) of $\cb_{i+1}$.  
Recall from \cite{ot} that $X(\cb_{i+1})=
X(\cb'_{i+1})\times \C^{\times}$.  
From the K\"unneth formula, we obtain the 
following exact sequence of $B'_{i+1}$-modules:
\begin{equation}
\label{eq:kunneth}
\xymatrix{
0\ar[r]&B'_{i+1}\ar[r]&B_{i+1}\ar[r]&B'_{i+1}(-1)\ar[r]&0.
}
\end{equation}

Let $J_{i+1}$ denote the kernel of the canonical projection 
$B_{i+1}\twoheadrightarrow A_{i+1}$. 
If we let $J'=J_{i+1}\cap B'_{i+1}$, 
then $J_{i+1}=B_{i+1}\otimes_{B'_{i+1}}J'$, as a 
module over $B_{i+1}$.   
Since $\ca$, $\cb$ are obtained from $\ca_{i+1}$, 
$\cb_{i+1}$, respectively,
by a sequence of fibred extensions, $J=B\otimes_{B_{i+1}} J_{i+1}$.

On the other hand, $B_{i+1}$ is a free module over $B'_{i+1}$, and
by applying Lemma~\ref{lem:free} inductively, $B$ is free over 
$B_{i+1}$.  Therefore, $B'_{i+1}\rightarrow B$ is a flat change
of rings, and it is enough to check that
\begin{equation}
\label{eq:extbjk}
\ext^p_{B'_{i+1}}(J',\k)_q=0
\end{equation}
if $q>d$.  By Lemma~\ref{lem:ExtAk}, $\ext^p_{B'_{i+1}}(J',\k)_q=
\ext^{p+1}_{B'_{i+1}}(A'_{i+1},\k)_{q-1}$.  Since $B'_{i+1}$ is Koszul
and $(A'_{i+1})_q=0$ for $q>d-1$, the rank of the arrangement, the 
groups \eqref{eq:extbjk} are zero for $q>d$ by \cite[Lemma 2.2]{joerg99}.
\end{proof}

The Lemma says, in particular, that the $B$-module $J(-c)$ 
has Castelnuovo-Mumford regularity no greater than the length 
of the singular range, $d-c$.  Moreover, the Lemma  
gives a combinatorial condition for the hypotheses of 
Theorem \ref{th:main} to be satisfied.  

\begin{corollary}
\label{cor:apply sing}
If $\ca$ is hypersolvable and its singular range has length $0$ 
or $1$, then $\g_A\cong \Lie(M[-2])\rtimes \h_A$.
\end{corollary}

\begin{example}[$2$-generic arrangements of rank $4$]
\label{ex:2gen4}
Suppose $\ca$ is a central arrangement in $\C^4$, with the property
that no three hyperplanes contain a common plane.  Such an 
arrangement is hypersolvable, by Theorem~\ref{th:slice}, with 
supersolvable deformation $\cb$ a Boolean arrangement.  
From Definition~\ref{def:srange} and Lemma~\ref{lem:cd}, 
$3\leq c\leq d\leq 4$, so the singular range has length 
$0$ or $1$.
\end{example}

On the other hand, the arrangement from 
Example~\ref{ex:lived2} is hypersolvable, with singular 
range $(3,5)$, and Corollary \ref{cor:apply sing} does not 
apply (indeed, its conclusion fails).

\subsection{Generic slices of supersolvable arrangements}
\label{subsec:hyper slice}

Lemma~\ref{lem:CMreg} provides bounds on the 
polynomial degrees of the homotopy module $M$, which 
cannot be improved without imposing further restrictions
on the arrangement.  In general, it is not obvious how 
to characterize the support of $M$ combinatorially; the 
problem seems similar to that of characterizing which 
arrangements have quadratic defining ideals, investigated 
in particular in~\cite{falk02,denyuz02}.  To this end, we 
isolate a class of hypersolvable arrangements for 
which the situation is more manageable. 

\begin{definition}
\label{def:gen slice}
A codimension-$k$ linear space $W$ is said to be 
generic with respect to an arrangement $\cb$ if 
$\rk(X\cap W)=\rk X+k$ for all $X\in L(\cb)$
with $\rk X\leq\rk \cb-k$. 
\end{definition} 

If $\cb$ is an essential,
supersolvable arrangement of rank $m$ and $W$ 
is a proper, linear space of dimension $\ell\geq3$, then by 
Theorem~\ref{th:slice}, the arrangement $\ca=\cb\cap W$ 
is hypersolvable.  We call such an arrangement a 
{\em generic (hypersolvable) slice} of rank $\ell$.

Not every hypersolvable arrangement is a generic slice, 
see Example 4.15 from \cite{ps02}.  

\begin{lemma}
\label{lem:srslice}
Let $\cb$ be a rank $m$ supersolvable arrangement, 
and let $\ca$ be a rank~$\ell$ generic slice. 
Then the singular range of $\ca$ is $(\ell,\ell)$.
\end{lemma}

\begin{proof}
The assumption of genericity means 
$L(\ca)_{\leq\ell-1}\cong L(\cb)_{\leq\ell-1}$.
However, $X\cap W=0$ for all $X\in L(\cb)_\ell$, so since 
$W$ is proper and $\cb$ is essential, the singular range of 
$\ca$ is $(\ell,d)$ for some $d$.  On the other hand, 
$\rk\ca_{\ell}=\rk\ca_m=\ell$, so the last $m-\ell$ 
extensions are all singular, and $d=\ell$.
\end{proof}

This is to say that, for generic slice arrangements, 
the module $J(-\ell)$ has a linear resolution.  Slightly 
more generally:

\begin{prop}
\label{prop:slice}
Let $\ca$ be a rank $\ell$ hypersolvable arrangement.  
Suppose there exists a generic slice $\mathcal{C}$ 
and fibred extensions
$\mathcal{C}=\ca_{m-i}\subset\cdots\subset\ca_{m-1}
\subset\ca_m=\ca$, for some $i\geq0$.  Then the singular 
range of $\ca$ is $(\ell,\ell)$.
\end{prop}

\begin{proof}
As in the proof of Lemma~\ref{lem:CMreg}, we may reduce 
to the case where $\ca=\mathcal{C}$, a generic slice of 
rank $\ell$.  Let $\cb$ be the supersolvable deformation of 
$\ca$.  Denote by $A'$ and $B'$ the Orlik-Solomon algebras of 
the respective decones, and let $J'=\ker (B'\twoheadrightarrow A')$.  

Let $R'=(B')^!$, and let $\ck=R'\otimes_\k (B')^*$ 
be the corresponding Koszul complex.  That is,
$\ck^q=R'(-q)\otimes_\k({B'}^{q})^*$ for $q\geq0$, 
with differential 
$\partial \colon e^*_i\otimes 1\mapsto 1\otimes x_i$. 
Since $B'$ is a Koszul algebra, $\ck$ is a free resolution 
of $\k$ over $R'$. 

Let $M'$ be the $\ell$-th syzygy module in the resolution 
$\ck\to  \k\to 0$.  That is, $M'$ is the cokernel of $\partial_{\ell+1}$, 
a left $R'$-module, which means $M'$ has minimal free resolution
\begin{equation}
\label{eq:resm}
\xymatrixcolsep{22pt}
\xymatrix{
0\ar[r] & \ck^m\ar[r]^{\partial_m} & \cdots\ar[r] & 
\ck^{\ell+1}\ar[r]^(.53){\partial_{\ell+1}}&\ck^\ell
\ar[r]^(.45){\eta} & M'\ar[r] & 0.
}
\end{equation}
From this we see that $M'$ is concentrated in internal
degree $\ell$, and $\ext_{R'}(M',\k)\cong J'$, as a $B'$-module.  
Since Koszul duality is an involution, $\ext_{B'}(J',\k)\cong M'$ 
as a left $R'$-module, and $M'$ is bigraded as claimed.
\end{proof}

The Proposition gives another criterion for the hypotheses
of Theorem \ref{th:main} to be satisfied.  We obtain:

\begin{corollary}
\label{cor:apply generic}
If $\ca$ is obtained by fibred extensions 
of a generic slice of a supersolvable arrangement, 
then $\g_A\cong \Lie(M[-2])\rtimes \h_A$.
\end{corollary}

\subsection{Hilbert series}
\label{subsec:hilb}
Expressions for the Hilbert series of the graded module 
$M=\ext_B(J,\k)$ are not known in general:  
compare with \cite{SS}.  However, a simple formula 
exists for generic slices, which can be extended to 
fibred extensions of generic slices.

Let $\beta_i$ denote the $i$th Betti number of $B'$, 
so that $h(B',t)=\sum_{i=0}^m\beta_i t^i$ is its Hilbert 
series.  The following fact is well-known; see~\cite{ot}.

\begin{lemma}
\label{lem:hsB}
There exist positive integers $1=d_1\leq d_2\leq\cdots\leq d_m$ 
for which
\begin{equation*}
h(B',t)=\prod_{j=2}^m(1+d_j t).
\end{equation*}
\end{lemma}

By taking the Euler characteristic of \eqref{eq:resm}, 
we note that for a generic slice of dimension $\ell$,
\begin{equation}
\label{eq:hilb M}
h_R(M,t)=h_{R'}(M',t)=
h(R',t)\sum_{i=0}^{m-\ell}(-1)^i \beta_{i+\ell} t^i.
\end{equation}
More generally, a fibred extension results in the same formula.
Under the hypotheses of Theorem~\ref{th:main}, together with
formula \eqref{eq:hilb M}, we have:

\begin{corollary}
\label{cor:hs}
If $h(U,t,u)=\sum_{p,q} \dim_\k U^p_qt^pu^q$ 
is the bigraded Hilbert series of $U=U(\g_A)$, then
\begin{equation}
\label{eq:hs1}
h(U,t,u)=h(R,t)\left(1-u^{-2}h_R(M,t,u)\right)^{-1}.
\end{equation}
In the case of a generic slice of dimension $\ell$,
\begin{equation}
\label{eq:hilb U}
h(U,t,u)=h(R,t)
\bigg(1-t^2u^{-2} h(R,t)
\sum_{i=0}^{m-\ell}(-1)^i\beta_{i+\ell}t^i\bigg)^{-1}.
\end{equation}
\end{corollary}

\section{A presentation for the homotopy Lie algebra}
\label{sec:lie}

For the hypersolvable arrangements satisfying the hypotheses of 
Theorem \ref{thm:main-arr}, the problem of writing an explicit 
presentation for the homotopy Lie algebra $\g_A$ 
is equivalent to that of presenting the homotopy 
module $M_A=\ext_B(J,\k)$.  We carry out this computation 
for generic slices of supersolvable arrangements.

Let $\ca=\{H_1,\dots, H_n\}$ be a hypersolvable arrangement, 
with supersolvable deformation $\cb$. 
As usual, let $\h$ denote the holonomy Lie algebra, and $R=U(\h)$ 
its enveloping algebra.  Recall $\h$ has a presentation with $n$ 
generators $x_1,\ldots, x_n$ in degree $(1,0)$, one for each 
hyperplane $H_i\in \ca$, and for each flat $F\in L_2(\ca)=L_2(\cb)$, 
relations
\begin{equation}
\label{eq:hol}
[x_i,\sum_{j\in F} x_j]=0,
\end{equation}
for all $i$ for which $i\in F$ (i.e., $F\subset H_i$).

Now assume $\ca$ is a generic slice of a supersolvable 
arrangement. Then the resolution \eqref{eq:resm} gives a 
presentation of the (deconed) homotopy module $M'$ 
as an $R'$-module.  In order to use this presentation 
explicitly, we will choose the basis for $B'^*$ given by 
identifying it with the flag complex of $\cb'$, for which 
we refer to \cite{denyuz02}.

Recall $\fl_p$ is a free $\k$-module on ``flags'' 
$(F_1,\ldots,F_p)$, where $F_i\in L_i(\cb')$ for 
$1\leq i\leq p$, and $F_i<F_{i+1}$, modulo the 
following relations:
\begin{equation}
\label{eq:flags}
\sum_{G: \, F_{i-1}<G<F_{i+1}}
(F_1,\ldots,F_{i-1},G,F_{i+1},\ldots,F_p),
\end{equation}
for each $i$, $1<i<p$.  Moreover, 
the map $f\colon \fl_p \to (B'^p)^*$ given by
\begin{equation}
\label{eq:fmap}
f \colon \left(F_1,\ldots,F_p\right) \mapsto 
\Big(\sum_{i\in F_1}e^*_i \Big)
\Big(\sum_{i\in F_2-F_1}e^*_i \Big)\cdots 
\Big(\sum_{i\in F_p-F_{p-1}}e^*_i \Big),
\end{equation}
is an isomorphism, cf.~\cite[dual of (2.3.2)]{sv}.

Under the identification $\fl\cong B'^*$, the boundary map 
in the Koszul complex becomes the following.  Given a flag 
$\F=(F_1,\ldots,F_p)$ and $i\in F_p$, define an element 
$\F-i\in \fl^{p-1}$ by finding the integer $j$ for which 
$i\in F_j-F_{j-1}$, and letting
\begin{equation}
\label{eq:Fi}
\F-i=(-1)^{j-1}\sum(F_1,\ldots,F_{j-1},G_j,G_{j+1},\ldots,G_{p-1}),
\end{equation}
where the sum is taken over all flags with the property that
$i\not\in G_{p-1}$ and $G_k<F_{k+1}$ for all $k$, $j\leq k<p$.  
Then the boundary map is given by extending
\begin{equation}
\label{eq:bdry map}
\partial \colon 
(F_1,\ldots,F_p)\mapsto\sum_{i\in F_p}(\F-i)\otimes x_i
\end{equation}
$R$-linearly.

For each element $\F\in \fl^{\ell}$, let $y_\F$ denote the 
corresponding element of $M'$; that is,
$y_\F=\eta\circ(f\otimes1)(\F\otimes1)$.  In particular,
we find a minimal generating set for $M'$ by choosing 
a set of $\beta_{\ell}$ flags of length $\ell$ 
in $L(\cb)$ appropriately.  In particular, one may construct a basis for 
$\fl^{\ell}$ using {\bf nbc}-sets: see, for example,
\cite[Lemma~3.2]{denyuz02}.

Then the relations in $M'$ are given by the image of 
$\partial_{\ell+1}$ in \eqref{eq:resm}.  We have, for 
each flag $\F=(F_1,\ldots,F_{\ell+1})$, a relation in $M'$ 
of the form
\begin{equation}
\label{eq:sumYX}
\sum_{i\in F_{\ell+1}} y_{\F-i}x_i.
\end{equation}
It follows that in $\g_A$, for each flag 
$\F=(F_1,\ldots, F_{\ell+1})$, we have a relation
\begin{equation}
\label{eq:rels}
\sum_{i\in F_{\ell+1}} [x_i,y_{\F-i}].
\end{equation} 

Now $M'$ is the restriction of the module $M$ from $R$ to $R'$, so the 
above gives a presentation for $M$ as well, noting that the central
element $\sum_{i=1}^n x_i$ in $R$ acts trivially.
One can find a minimal set of relations just by taking the flags 
$\F$ above to come from a basis of $\fl^{\ell+1}$.  
We summarize this discussion, as follows. 

\begin{theorem}
\label{thm:g pres}
Let $\ca=\{H_1,\dots, H_n\}$ be a generic slice of 
a supersolvable arrangement, and let $A$ be the 
Orlik-Solomon algebra of $\ca$.  Then, the homotopy 
Lie algebra $\g_A$ has presentation with generators 
\begin{itemize}
\item 
$x_i$ in degree $(1,0)$, for each $i\in [n]$,
\item
$y_\F$ in degree $(2,\ell-2)$, for each $\F\in \fl^\ell$,
\end{itemize}
and relations 
\begin{itemize}
\item 
$\big[ x_i,\sum_{j\in F} x_j \big] =0$, for each flat $F\in L_2(\ca)$ 
and each $i\in F$, 
\item 
$\sum_{i\in F_{\ell+1}}\big[ x_i,y_{\F-i}\big] =0$, 
for each flag $\F=(F_1,\ldots, F_{\ell+2})\in \fl^{\ell+1}$,
\item
$\big[\sum_{i=1}^n x_i,y_\F\big]=0$, for each $\F\in\fl^\ell$.
\end{itemize}
\end{theorem}

We illustrate the above with an example.
\begin{example}
\label{ex:pres}
Consider the arrangement $\ca$ defined by the polynomial
\[
Q_\ca=xyz(x-z)(y-z)(2x-y-4z)(2x-y-5z)(x+5y+2z)(x+5y+z).
\]
This is a generic slice of the supersolvable arrangement $\cb$, 
the cone over the arrangement defined by the polynomial
$Q_{\cb'}=vwxy(x-1)(y-1)(v-1)(w-1)$.  
The Poincar\'e polynomials of the deconed arrangements 
are given by 
\begin{align*}
\pi(\ca',t)&=1+8t+24t^2, \\ 
\pi(\cb',t)&=(1+2t)^4=1+8t+24t^2+32t^3+ 16t^4.
\end{align*}
Thus the homotopy module $M'$ has $32$ generators
and $16$ relations, which can be described as follows.

Label the hyperplanes of $\cb'$ as 
$0_0,1_0,2_0,3_0,0_1,1_1,2_1,3_1$, in the order above.  
A basis of $32$ flags of length $3$ can be constructed by choosing
three intersecting hyperplanes $i_a,j_b,k_c$, with $0\leq i<j<k\leq 3$
and $a,b,c\in\set{0,1}$, and forming a flag by successively intersecting 
the hyperplanes, from right to left.  We will call this flag $\F_{i_aj_bk_c}$.
A basis of $16$ flags of length $4$ in $\cb'$ is constructed by choosing
four intersecting hyperplanes, $0_a,1_b,2_c,3_d$, for all choices of
$a,b,c,d\in\set{0,1}$, and forming a flag again by successive intersection.

Let $\g_A$ be the holonomy Lie algebra of $\ca$.  
Then $\g_A$ has one generator $x_H$ for each hyperplane $H$,
together with $32$ additional generators $y_{i_aj_bk_c}$ in 
degree $(2,1)$, and relations
\[
[x_{0_a},y_{1_b2_c3_d}]-[x_{1_b},y_{0_a2_c3_d}]+[x_{2_c},y_{0_a1_b3_d}]-
[x_{3_d},y_{0_a1_b2_c}],
\]
for each $a,b,c,d\in\set{0,1}$, in addition to the holonomy relations
\eqref{eq:hlie arr}, and relations
\[
\Big[\sum_{H\in\ca}x_H, y_{i_aj_bk_c} \Big]
\]
for each choice of $i,j,k$, $a,b,c$.
\end{example}

\section{Two-generic arrangements of rank four}
\label{sec:two generic}

We now present a method for computing the Hilbert series 
of the homotopy Lie algebra of a particularly nice class 
of arrangements:  rank-$4$ arrangements for which no 
three hyperplanes contain a common plane.  

For any rank $\ell$ arrangement $\ca$ with $n$ hyperplanes, 
let $E=\bigwedge_{\k} (e_1,\dots, e_n)$ be the 
exterior algebra, $A=E/I$ the Orlik-Solomon algebra, 
and $S=\k[x_1,\ldots,x_n]$ the polynomial algebra. 
We recall the following.  

\begin{theorem}[Eisenbud-Popescu-Yuzvinsky \cite{epy03}]
\label{thm:epy}
The complex of $S$-modules
\[
\xymatrix{
0& F(\ca)\ar[l] & A^\ell\otimes S\ar[l] & \ldots\ar[l] &
A^1\otimes S\ar[l] & A^0\otimes S\ar[l] & 0\ar[l]
}
\]
is exact, where the boundary maps are induced by 
multiplication by $\sum_{i=1}^n e_i\otimes x_i$, and 
the $S$-module $F(\ca)$ is taken as the cokernel of 
the map $A^{\ell-1}\otimes S\to  A^\ell\otimes S$.
\end{theorem}

It follows from Bernstein-Gelfand-Gelfand duality that, 
for each $p\geq0$, there is a graded isomorphism of 
$S$-modules,
\begin{equation}
\label{eq:ext fas}
\ext_E^p(A,\k)_q=\ext_S^{\ell-q}(F(\ca),S)_{p+q}.
\end{equation}
We refer to \cite{SS} for the case of the smallest $q>0$
for which this is nonzero.  Details will appear in further work.

Now let $\ca$ be a $2$-generic arrangement.
Notice that $B=E$ and $U(\h)=B^!=S$.  Then, applying Lemma~\ref{lem:ExtAk}
to \eqref{eq:ext fas}, we obtain
\begin{equation}
\label{eq:bgg}
M^p_q=\ext_S^{\ell-q+1}(F(\ca),S)_{p+q},
\end{equation}
for $p\geq0$ and $0\leq q\leq\ell$.
As a result, presentations for the $S$-modules $M_q$ can 
be obtained computationally for specific examples, using 
formula \eqref{eq:bgg}. 

We recall from Example~\ref{ex:2gen4} that, if the rank of the arrangement
$\ell=4$, then $\ca$
satisfies hypotheses \eqref{i} of Theorem~\ref{thm:main-arr}:
$M_q=0$ unless $q=3$ or $q=4$, i.e., the singular range of $\ca$ 
is $(3,4)$.  

\begin{example}
\label{ex:2gen-7}

Consider arrangements $\ca_1$ and $\ca_2$ defined by the 
polynomials
\begin{align*}
Q_1&=xyzw(x+y+z)(y+z+w)(x-y+z+w), \\
Q_2&=xyzw(x+y+z)(y+z+w)(x-y+z-w).
\end{align*}
Both arrangements have $7$ hyperplanes and $5$ lines that 
each contain $4$ hyperplanes, so the characteristic polynomials 
are $\pi(\ca_1,t)=\pi(\ca_2,t)=1+7t+21t^2+30t^3+15t^4$.  
Since there are no nontrivial intersections in codimension $2$, 
the fundamental group of both complements is $\Z^7$,
and $R=U(\h)$ is a polynomial ring.  

We now use \eqref{eq:bgg} to compute the Hilbert series of the 
graded modules $M_3$ and $M_4$ (recalling $M_q=0$ for $q\neq3,4$). 
With the help of Macaulay 2, we find for $\ca_1$ 
\begin{align*}
h(M_3,t)&=(5+ 2t)/(1-t)^3 
=5 + 17t + 36 t^2 + 62 t^3 +\cdots \\
h(M_4,t)&=(2 - t)(1 + 2t + 2t^2)/(1 - t)^6 
= 2 + 15 t+ 62 t^2+ 185 t^3 +\cdots  \\
\intertext{while for $\ca_2$,}
h(M_3,t)&=(5+ t)/(1-t)^3 
=5+ 16t+33t^2+ 56t^3 +\cdots \\
h(M_4,t)&=(1 + 6t - t^2 - t^3)/(1 - t)^6
=1+12t+ 56t^2+ 175t^3 +\cdots 
\end{align*}

Using formula \eqref{eq:hs1}, this yields 
expressions for the Hilbert series of $U(\g_1)$ and $U(\g_2)$.  
Comparing these Hilbert series shows $U(\g_1)\not\cong U(\g_2)$, 
and hence the two arrangements must have 
non-isomorphic homotopy Lie algebras.
\end{example}

\begin{example}
\label{ex:blocks}
In 1946, Nandi~\cite{nandi46} showed that there are 
exactly three inequivalent block designs with parameters 
$(10,15,6,4,2)$.  We list the blocks of each below.  
Each block design gives rise to a rank-$4$ matroid
on ten points by taking the dependent sets to be 
those subsets that either contain one of the blocks 
or contain at least five elements.

\[
\begin{array}{|l|l|} \hline
D_1 & \{{\tt abcd,abef,aceg,adhi,bchi,bdgj,cdfj,afhj,agij,}\\
& {\tt behj,bfgi,ceij,cfgh,defi,degh}\} \\ \hline
D_2 & \{{\tt abcd,abef,aceg,adhi,bcij,bdgh,cdfj,afhj,agij,}\\
& {\tt aehj,bfgi, cehi, cfgh, defi, degj}\} \\ \hline
D_3 & \{{\tt abcd,abef,acgh,adij,bcij,bdgh,cdef,aegi,afhj,}\\
& {\tt behj,bfgi,cehi,cfgj,degj,dfhi}\} \\ \hline 
\end{array}
\]
\vskip 8pt
By construction, there are no nontrivial, dependent sets of size three,
so each arrangement is $2$-generic.

If we call the corresponding Orlik-Solomon algebras $A_1$, $A_2$, 
and $A_3$, it is straightforward to calculate that 
$h(A_i,t)=1+10t+45t^2+105t^3+69t^4$ for $i=1,2,3$.  
In each case, the singular range is $(3,4)$.  The ideals
$J_1$, $J_2$, $J_3$ have differing resolutions, however, from which
it follows that $\g_{A_1}$, $\g_{A_2}$, and $\g_{A_3}$ are pairwise
non-isomorphic.
\end{example}

\section{Topological interpretations}
\label{sect:top}

\subsection{Generic slices}
\label{subsec:generic slices}
A particularly simple situation, analyzed in detail by 
Dimca and Papadima in \cite{dimpap03}, 
is when $\ca$ is a generic slice of rank $\ell>2$ of a 
supersolvable arrangement $\cb$.  Let $\ca'$ and $\cb'$ 
be the respective decones, with complements 
 $X=X(\ca')$ and $Y=X(\cb')$.  The two 
spaces share the same fundamental group, $\pi$, 
and the same integral holonomy Lie algebra, $\h$. 

In \cite[Theorems 18(ii) and 23]{dimpap03}, Dimca 
and Papadima establish the following facts.   
The universal enveloping algebra $U(\h)$ 
is isomorphic (as a Hopf algebra) to the 
associated graded algebra $\gr_{I\pi}(\Z\pi)$, where 
$\Z\pi$ is the group ring of  $\pi$, 
with filtration determined by the powers 
of the augmentation ideal $I\pi$.   
The first non-vanishing higher homotopy group of $X$ 
is $\pi_{\ell-1}(X)$; when viewed as a module over $\Z\pi$, 
it has resolution of the form
\begin{equation}
\label{pipres}
\xymatrixcolsep{16pt}
\xymatrix{
0\ar[r] & H_m(Y)\otimes\Z\pi \ar[r]&
\cdots \ar[r]&H_{\ell}(Y)\otimes\Z\pi\ar[r]&
\pi_{\ell -1}(X)\ar[r]& 0
}.  
\end{equation}
Finally, the associated graded module of $\pi_{\ell-1}(X)$, 
with respect to the filtration by powers of $I\pi$, has Hilbert series
\begin{equation}
\label{hilbpip}
h(\gr_{I\pi}^\bfcdot\pi_{\ell-1}(X),t)=
(-1/t)^{\ell}\left(1-\frac{\sum_{j=0}^{\ell-1} 
(-1)^j\beta_jt^j}{%
\sum_{j=0}^m(-1)^j\beta_jt^j}\right), 
\end{equation}
where $\beta_j$ are the Betti numbers of $Y$. 

Consider the integral cohomology rings 
$A=H^\bfcdot(X,\Z)$ and $B=H^\bfcdot(Y,\Z)$.  
We have $(B^i)^*=H_i(Y,\Z)$, since the homology of an arrangement 
complement is torsion-free.  Thus, tensoring with $\k$, and passing 
to the associated graded in resolution \eqref{pipres} recovers 
resolution  \eqref{eq:resm}.  As a consequence, we obtain the 
following.

\begin{prop}
\label{prop:homotopy}
Let $\ca$ be a generic slice of rank $\ell>2$ of a 
supersolvable arrangement, and let  $X=X(\ca')$ be 
the complement of its decone.  
The homotopy module of the algebra $A=H^\bfcdot(X,\k)$ 
is isomorphic to the graded module associated to the the 
first nonvanishing higher homotopy group of $X$:
\begin{equation}
\label{mpip}
M_A \cong \gr_I^\bfcdot\pi_{\ell-1}(X) \otimes \k.
\end{equation}
\end{prop}

\subsection{Rescaling}
\label{subsec:rescaling}
Fix an integer $q\ge 1$.  
The $q$-rescaling of a graded algebra $A$ 
is the graded algebra $A^{[q]}$, with 
$A^{[q]}_{i(2q+1)} =A_i$ and $A^{[q]}_{j}=0$ if $(2q+1)\nmid j$,  
and with multiplication rescaled accordingly. When taking 
the Yoneda algebra of $A^{[q]}$, the internal degree of 
the Yoneda algebra of $A$ gets rescaled, 
while the resolution degree stays unchanged:
\begin{equation}
\label{eq:resc yoneda}
\ext_{A^{[q]}}(\k,\k) = \ext_{A}(\k,\k)^{[q]}.
\end{equation}

Similarly, the $q$-rescaling of a graded Lie algebra $L$ 
is the graded Lie algebra $L^{[q]}$, with  $L^{[q]}_{2iq} =L_i$ 
and $L^{[q]}_j=0$ if $2q\nmid j$, and with Lie bracket rescaled 
accordingly.  Rescaling works well with the holonomy 
and homotopy Lie algebras: 
\begin{equation}
\label{eq:resc hg}
\h_{A^{[q]}} = \h_A^{[q]}, \quad 
\g_{A^{[q]}} = \g_A^{[q]}.
\end{equation}
The Hilbert series of the enveloping algebras of $\g_A^{[q]}$ 
and $\g_A$ are related as follows:
\begin{equation}
\label{eq:hrescale}
h(U(\g_A^{[q]}),t,u)=h(U(\g_A),tu^{2q},u^{2q+1}).
\end{equation}

Now let $X$ be a connected, finite-type CW-complex.  
A simply-connected, finite-type CW-complex $Y$ is 
called a $q$-rescaling of $X$ (over a field $\k$) if the 
cohomology algebra $H^{\bfcdot}(Y,\k)$ is the 
$q$-rescaling of $H^{\bfcdot}(X,\k)$, i.e.,
\begin{equation}
\label{eq:resc}
H^{\bfcdot}(Y,\k)=H^{\bfcdot}(X,\k)^{[q]}.
\end{equation}

Rational rescalings always exist:  take a Sullivan 
minimal model for the $1$-connected, finite-type 
differential graded algebra $(H^{\bfcdot}(X,\Q)^{[q]}, d=0)$, 
and use \cite{Su} to realize it by a finite-type, 
$1$-connected CW-complex, $Y$.  
The space constructed this way is the desired 
rescaling. Moreover, $Y$ is formal, i.e, its rational
homotopy type is a formal consequence of its rational 
cohomology algebra.  Hence, $Y$ is uniquely determined, 
up to rational homotopy equivalence, among spaces 
with the same cohomology ring (though there may be other, 
non-formal rescalings of $X$, see \cite{ps04}).

\begin{prop}
\label{prop:resc}
Let $X$ be a finite-type CW-complex, with 
cohomology algebra $A=H^{\bfcdot}(X;\Q)$.  
Let $Y$ be a finite-type, simply-connected 
CW-complex with $H^{\bfcdot}(Y;\Q)\cong A^{[q]}$. 
If $Y$ is formal, then
\begin{equation}
\label{eq:holo pi}
\pi_{\bfcdot}(\Omega Y) \otimes \Q \cong \g_{A}^{[q]}. 
\end{equation}
\end{prop}

\begin{proof}
Since $Y$ is formal, the Eilenberg-Moore spectral 
sequence of the path fibration $\Omega Y \to P Y \to Y$ 
collapses, yielding an isomorphism of Hopf algebras 
between the Yoneda algebra of $H^{\bfcdot}(Y;\Q)$ 
and the Pontryagin algebra $H_{\bfcdot}(\Omega Y; \Q)$.   
From the rescaling assumption, we obtain
\begin{equation}
\label{eq:ext pont}
\ext_{A^{[q]}} (\Q,\Q) \cong H_{\bfcdot}(\Omega Y; \Q),
\end{equation}
By Milnor-Moore \cite{mimo65}, we find that 
$\g_{A^{[q]}} \cong \pi_{\bfcdot}(\Omega Y) \otimes \Q$, 
as Lie algebras.  Using \eqref{eq:resc hg} finishes 
the proof.
\end{proof}

As a consequence, we obtain a quick proof of a special 
case of Theorem A from \cite{ps04}.  

\begin{corollary}[\cite{ps04}]
\label{rem:both formal}
Suppose $X$ and $Y$ are spaces as above. 
If both $X$ and $Y$ are formal and $A$ is Koszul, then
\begin{equation}
\label{eq:resc pi}
\pi_{\bfcdot}(\Omega Y)\otimes \Q \cong
\left(\gr_{\bfcdot}(\pi_1 X)\otimes \Q\right)^{[q]}. 
\end{equation}
\end{corollary}

\begin{proof}
Since $A$ is Koszul, $\g_{A}=\h_{A}$. Since $X$ is formal, 
$\gr_{\bfcdot}(\pi_1 X)\otimes \Q\cong \h_{A}$, cf.~\cite{Su}. 
The conclusion follows from \eqref{eq:holo pi}.
\end{proof}

\begin{remark}
\label{rem:rht}
When $X$ is formal (but not necessarily simply connected), 
a theorem of Papadima and Yuzvinsky~\cite{py} states that the 
cohomology algebra $A=H^\bfcdot(X;\Q)$ is Koszul if 
and only if the Bousfield-Kan rationalization $X_\Q$ is 
aspherical. Now, by a classical result of Quillen \cite{Qu68},  
$U(\h_A) \cong \gr_{I\pi}\Q\pi_1(X_\Q)$. More generally, 
it seems likely that
\begin{equation}
\label{eq:ft}
U(\g_A)\cong U(\pi_{\bfcdot}(\Omega \widetilde{X}_\Q))
\,\widehat{\otimes}\, \gr_{I\pi}\Q\pi_1(X_\Q),
\end{equation}
in view of a result of F\'elix and Thomas~\cite{feltho86}.  
(Here again, $\Q\pi_1(X_\Q)$ acts on the left-hand factor 
by the action induced from $\pi_1(X_\Q)$ on the universal 
cover $\widetilde{X_\Q}$.)

However, if $X$ is a hyperplane arrangement 
complement, then $X$ is not in general a nilpotent space.   
This means that we can expect to find such spaces $X$ for 
which $\pi_i(X_\Q)\not\cong\pi_i(X)\otimes\Q$.  The first such 
example was found by Falk \cite{falk88}, who noted that the 
complement $X$ of the $D_4$ reflection arrangement 
is aspherical, while its Bousfield-Kan rationalization
$X_\Q$ is not.  In general, then, we know of no way to 
relate $\g_A$ with the topological homotopy Lie algebra, 
$\pi_{\bfcdot}(\Omega X)\otimes \Q$.
\end{remark}

\subsection{Redundant subspace arrangements}
\label{subsec:redundant}
Let $\ca=\{H_1,\dots,H_n\}$ be an arrangement of hyperplanes 
in $\C^\ell$.  If $q$ is a positive integer, then 
$\ca^{(q)}=\{H_1^{\times q},\dots,H_n^{\times q}\}$ 
is an arrangement of codimension $q$ subspaces in $\C^{q\ell}$. 
For example, if $\ca$ is the braid arrangement in $\C^{\ell}$,
with complement equal to the configuration space
of $\ell$ distinct points in $\C$, then the complement 
of $\ca^{(q)}$ is the configuration space of $\ell$ distinct 
points in $\C^q$. 

\begin{prop}
\label{prop:red resc}
Let $\ca$ be a hyperplane arrangement, with Orlik-Solomon 
algebra $A=H^{\bfcdot}(X;\Q)$.  Fix $q\ge 1$, and let 
$Y=X(\ca^{(q+1)})$ be the complement of the corresponding 
subspace arrangement.  Then:
\begin{equation}
\label{eq:resc hyper}
\pi_{\bfcdot}(\Omega Y) \otimes \Q \cong \g_{A}^{[q]}.
\end{equation}
\end{prop}

\begin{proof}
Clearly, $Y$ is simply-connected.  As shown in \cite{CCX}, 
$H^{\bfcdot}(Y;\Q)$ is the $q$-rescaling of $H^{\bfcdot}(X;\Q)$. 
Since $\ca^{(q+1)}$ has geometric intersection 
lattice, its complement $Y$ is formal, see 
\cite[Prop.~7.2]{yuz02}. The conclusion then 
follows from Proposition~\ref{prop:resc}.
\end{proof}

\begin{corollary}
\label{cor:hyp resc}
Let $\ca$ be a hypersolvable arrangement, satisfying either 
of the hypothesis of Theorem \ref{thm:main-arr}.  Then 
\[
\pi_{\bfcdot}(\Omega Y) \otimes \Q \cong (\Lie(M_A[-2])\rtimes \h_A)^{[q]}.
\]
\end{corollary}

\begin{example}
\label{ex:red}
Let $\ca_1$ and $\ca_2$ be the hyperplane arrangements 
from Example~\ref{ex:2gen-7}.  Denote by $\g_i=\g_{A_i}$ 
the respective homotopy Lie algebras, $i=1,2$.  
Consider the redundant subspace arrangements 
$\ca_1^{(2)}$ and $\ca_2^{(2)}$.  Both are arrangements 
of $7$ codimension-$2$ complex subspaces of $\C^8$.  
Denoting their complements by $Y_1$ and $Y_2$, respectively, 
we have $\pi_1(Y_1)=\pi_1(Y_2)=0$ and $H_*(Y_1)\cong H_*(Y_2)$ 
as graded abelian groups.

Let $\pi_{\bfcdot}(\Omega Y_i)\otimes\Q$ be the respective 
(topological) homotopy Lie algebras. By Proposition~\ref{prop:resc}, 
we have $\pi_{\bfcdot}(\Omega Y_i)\otimes\Q\cong \g^{[1]}_{i}$.  
Making use of the previous calculations for the arrangements 
$\ca_1$ and $\ca_2$, together with formula \eqref{eq:hrescale}, 
we find that $U(\g_i^{[1]})_p$ has rank $1,0,7,0,28,0,84,5,210$ 
for $0\leq p\leq 8$, for both $i=1,2$. 
It follows that, for $p\leq9$, the group $\pi_p(Y_i)\otimes\Q=0$, 
except for $\pi_3(Y_i)\otimes\Q\cong\Q^7$ and 
$\pi_8(Y_i)\otimes\Q\cong\Q^5$. 

However, for $p=9$, the ranks of  $U(\g_i^{[1]})_p$ 
are $52$ and $51$, respectively. Hence, 
\[
\pi_{10}(Y_1)\otimes\Q\cong\Q^{17} \quad \text{and} \quad 
\pi_{10}(Y_2)\otimes\Q\cong\Q^{16}, 
\]
and so $Y_1\not\simeq Y_2$. 
\end{example}

\begin{ack}
We thank Srikanth Iyengar for helpful conversations.  
A substantial portion of this work was carried out while 
the authors were attending the program ``Hyperplane 
Arrangements and Applications'' at the Mathematical 
Sciences Research Institute in Berkeley, California, 
in Fall, 2004.  We thank MSRI for its support and 
hospitality during this stay. 
\end{ack}

\providecommand{\bysame}{\leavevmode\hbox to3em{\hrulefill}\thinspace}
\providecommand{\MR}{\relax\ifhmode\unskip\space\fi MR }
\providecommand{\MRhref}[2]{%
  \href{http://www.ams.org/mathscinet-getitem?mr=#1}{#2}
}
\providecommand{\href}[2]{#2}

\end{document}